# Multiwell rigidity in nonlinear elasticity

Milena Chermisi and Sergio Conti

*Fachbereich Mathematik, Universität Duisburg-Essen,*
*Lotharstrasse 65, 47057 Duisburg, Germany*

February 2, 2008

## 1 Introduction

Variational models from nonlinear elasticity have the form

$$\int_\Omega W(\nabla u)dx, \qquad (1.1)$$

where the energy density $W$ is invariant under rotations; in the simplest cases $W$ is minimized by the set of proper rotations $SO(n)$. The set $SO(n)$, much as the set of conformal matrices, is *rigid*, in the sense that there are no nontrivial gradient fields taking values in $SO(n)$. This classical result, due for smooth functions to Liouville, permits to show that minimizers of (1.1) are affine (on connected domains). Several improvements have been obtained, among others by Gehring [19], John [20], Reshetnyak [27], and Kohn [22], culminating in the recent quantitative version by Friesecke, James and Müller [18], who have shown that gradient fields taking values close to $SO(n)$ are approximately constant, in the sense that for any $u \in W^{1,2}(\Omega \subset \mathbb{R}^n; \mathbb{R}^n)$ one has

$$\min_{Q \in SO(n)} \int_\Omega |\nabla u - Q|^2 dx \le c \int_\Omega \mathrm{dist}^2(\nabla u, SO(n)) dx. \qquad (1.2)$$

This shows that low-energy states are approximately affine, with an estimate which has the optimal norm and scaling, and constitutes the nonlinear counterpart of the classical linear estimate known as Korn's inequality.

In the study of solid-solid phase transitions and in particular of shape-memory alloys one is interested in energy densities $W$ which are minimized by several copies of $SO(n)$, i.e., by sets of the form

$$K := SO(n)U_1 \cup \cdots \cup SO(n)U_l \qquad (1.3)$$

where $U_i \in \mathbb{R}^{n \times n}$ are the eigenstrains of the different phases, see, e.g., [2, 11, 3, 25, 13, 4, 26]. These sets are in general not rigid, and there are nontrivial gradient fields taking values in $K$, hence no estimate like (1.2) can be expected. Such gradient fields are however strongly restricted, and it is easy to see that if the interface between a region where $\nabla u = A$ and one where $\nabla u = B$ is sufficiently smooth, then $A - B = a \otimes \nu$, where $\nu$ is the normal to the interface. This implies that the interface is locally a hyperplane, and that its normal has to take one of finitely many directions, which are determined by $K$. More precisely, for the case of two wells in two dimensions, it was shown by Dolzmann and Müller [14] that if $\nabla u \in BV(\Omega; K)$ then $\nabla u$ is piecewise constant, and its jump set is the disjoint union of segments with endpoints on $\partial\Omega$, with only two possible orientations. Kirchheim was later able to obtain a corresponding result for the substantially more complex, but physically more relevant, three-well problem in three dimensions [21].

A first quantitative version of the Dolzmann-Müller rigidity result was obtained by Lorent [23] for the case of two wells with equal determinant, and bi-Lipschitz maps $u : B_1 \to \mathbb{R}^2$. Precisely, he has proven that if $\|D^2 u\|(B_1)$ is sufficiently small, then

$$\min_{J \in \{U_1, U_2\}} \|\text{dist}(\nabla u, SO(2)J)\|_{L^1(B_\rho)} \leq c \big(\|\text{dist}(\nabla u, SO(2)\{U_1, U_2\})\|_{L^1(B_1)}\big)^\gamma$$

for some constants $\gamma$, $\rho$, $c > 0$. Here it is tacitly assumed that $\nabla u \in BV$, and $D^2 u$ denotes the distributional second gradient. Conti and Schweizer [9] improved the estimate by showing that

$$\min_{J \in \{U_1, U_2\}} \|\text{dist}(\nabla u, SO(2)J)\|_{L^1(\Omega')} \leq c \|\text{dist}(\nabla u, SO(2)\{U_1, U_2\})\|_{L^1(\Omega)}, \qquad (1.4)$$

for any two matrices $U_1$ and $U_2$ with positive determinant and for any maps $u \in W^{1,1}(\Omega; \mathbb{R}^2)$ for which $\|D^2 u\|(\Omega)$ is small compared to $\text{dist}(\Omega', \partial \Omega)$, where $\Omega'$ is a connected subset of $\Omega$. A related estimate in a geometrically linear context had been obtained in [10]. These results have had several applications, e.g., in studying the scaling of singularly perturbed problem under Dirichlet boundary conditions [24, 5] or in proving compactness and $\Gamma$-convergence for a sequence of singularly perturbed functionals of the kind

$$I_\varepsilon[u] := \int \frac{1}{\varepsilon} W(\nabla u) + \varepsilon |\nabla^2 u|^2 dx,$$

see [9, 10]. In particular, rigidity estimates are needed in order to pass from the case where the nonconvex term of the energy has only finitely many minimizers [17, 7] to the elasticity case where the set $K$ of minimizers of $W$ is infinite.

We present here a generalization of (1.4) beyond the two-well, two-dimensional case. We shall assume that the wells are well separated, in the following sense:

$$\text{for each } i = 1 \ldots l \text{ there is } \xi_i \in \mathbb{S}^{n-1} \text{ such that } |U_i \xi_i| > \max_{j \neq i} |U_j \xi_i|. \qquad (1.5)$$

We shall show later that this condition holds in the most relevant examples, including the three-well problem in three dimensions, and any two-well problem with two rank-one connections (see Remarks 3.6, 3.7 and 3.8). Necessity of the assumption (1.5) is shown in Remark 3.9.

**Theorem 1.1.** *Let $U_1, \ldots U_l \in \mathbb{R}^{n \times n}$ have positive determinant, fulfill the separation condition (1.5), and let $K$ be defined as in (1.3). Let $\Omega' \subset\subset \Omega \subset \mathbb{R}^n$ be two bounded Lipschitz domains, with $\Omega'$ connected. Then there are positive constants $\eta_0$, $c_0$, $c_1$ and $c_2$ such that for any $u \in W^{1,1}(\Omega; \mathbb{R}^n)$ with $\nabla u \in BV$, satisfying*

$$\int_\Omega |D^2 u| \leq \eta_0 \, \text{dist}^{n-1}(\Omega', \partial \Omega) \qquad (1.6)$$

*where $\eta_0 = \eta_0(\{U_i\}, n)$, one has*

$$\min_{J \in \{U_i\}_i} \int_{\Omega'} \text{dist}(\nabla u, SO(n)J) \, dx \leq c_0 \int_\Omega \text{dist}(\nabla u, K) \, dx, \qquad (1.7)$$

*where $c_0 = c_0(\Omega, \Omega', \{U_i\}, n)$. Further, if*

$$\int_\Omega \text{dist}(\nabla u, K) \, dx \leq c_2(\{U_i\}, n) \text{dist}^n(\Omega', \partial \Omega), \qquad (1.8)$$

*then (1.7) holds with $c_0 = c_1(\{U_i\}, n)$, independently of the domains.*

The proof is separated into two parts: first one proves rigidity along segments for many segments (Section 2), and then one deduces $L^1$ rigidity from the rigidity of segments via geometrical arguments (Section 3). Whereas the subdivision of the argument in these two basic steps is the same as in the two-dimensional case [9], the proof of both parts presents significant differences and new difficulties with respect to the two-dimensional case.

The multiwell rigidity result on segments is presented in Proposition 2.1. At variance with [9], where the Sobolev embedding $W^{1,1}(\mathbb{S}^1) \subset L^\infty(\mathbb{S}^1)$ automatically gave uniform estimates on spheres, permitting a direct usage of Brouwer-degree theory to prove invertibility, here the same holds only if higher integrability is first obtained. This is achieved via a truncation argument, and by carefully passing from the original function to the truncated one and back in appropriate points of the proof. In particular, the map $u$ is replaced by a Lipschitz map $v$, such that $u = v$ outside a set of small perimeter. This requires a truncation argument which estimates also the perimeter, and not only the volume, of the set where the function is changed (Proposition 2.2).

The step from the segment rigidity of Proposition 2.1 to Theorem 1.1 is then essentially geometric. Whereas in two dimensions one could use two coupled triangles to obtain the optimal scaling in a region close to the midplane of the rhombus (see Figure 3 below), in the $n$ dimensional case we need to consider a small region close to the barycenter of a simplex. In both cases the key idea is that a direction $\xi_0$ can be chosen (according to (1.5)) so that the "majority" phase stretches that direction more than the others. Then, if one considers segments which have approximately that orientation inside a "rigid" simplex, the volume fraction of the "minority" phases can be controlled by the total stretch, i.e., by the multiwell energy. Details are discussed in Section 3.

In closing, we recall that a quite different situation arises in the case that the wells are incompatible. Indeed, in that case (1.2) can be directly generalized to the multiwell situation, as was shown by Chaudhuri and Müller [6]; a simpler proof of the same result was then obtained by De Lellis and Székelyhidi [12].

## 2 Segment rigidity

We derive in this section a segment rigidity result for multiwell energies and dimension $n \geq 2$, generalizing [9, Proposition 2.2]. Before giving the precise statement we briefly sketch the main ideas. Let $u : \Omega \to \mathbb{R}^n$ be such that $\nabla u$ is close, in $L^1$, to the set $K$ introduced above. Then the domain can be subdivided into $l$ parts, according to which well $\nabla u(x)$ is closer to. We assume that, in an appopriate weak formulation, one of those sets is large and the others have small perimeter (see (2.1) and (2.2)). Then a slicing argument shows that most segments do not intersect the "minority" sets, and have $\nabla u$ close to $SO(n)U_1$ (possibly after relabeling). Therefore the tangential component $\partial_\tau u$ of $u$ has length close to $|U_1 \tau|$, and hence the map $u$ does not "make the segment longer". In order to obtain the converse inequality, one applies a similar argument to the inverse of $u$.

The main difficulty is that $u$ is, in general, not invertible. This was overcome in [9, Proposition 2.2] by an argument based on Brouwer degree, exploiting an uniform estimate on the restriction of $u$ to the (one-dimensional!) boundary of the domain. The same estimate does not, however, hold in higher dimension, unless stronger assumptions on the integrability of $\nabla u$ are made. Such assumptions are however not expected to hold in the typical applications of this type of rigidity estimate discussed after (1.4).

To overcome this difficulty we replace the map $u$ by a suitable truncation $v$. The new map $v$ is Lipschitz, and agrees with $u$ away from a small set with small perimeter (see Proposition 2.2 for the truncation). Then we can obtain for $v$ the uniform estimate on the boundary, and obtain appropriate invertibility of $v$. Finally, it remains to show that most segments do not intersect the set where $u$ differs from $v$, and this is done using the fact that the latter set has small perimeter. Having been able to invert $u$ along the relevant segments, we can apply the "no-stretch" argument to the inverse and conclude that $u$ does not "make the segment shorter", which concludes the proof.

The function $u$ is assumed to be $C^1$. This permits in particular to use the implicit function theorem to prove local invertibility, and to have pointwise values for $\nabla u$. None of the constants entering the statement can depend on the $C^1$ norm of $u$. In particular the Lipschitz norm of the truncation $v$ does not depend on $u$, but only on $K$ and the constants $\bar{c}, p$ entering the statement.

**Proposition 2.1.** *Let $\alpha \in (0, 1/8)$, $\theta \in (0, 1)$, $\bar{c} > 0$ and $p \geq 1$. Let $K$ be defined as in (1.3) where $U_1, \ldots, U_l \in \mathbb{R}^{n \times n}$ have positive determinant. Then there are $\eta, c > 0$ depending only on the above quantities such that the following holds:*
*Let $\Omega \subset \mathbb{R}^n$ be open and convex, $x_0, y_0 \in \Omega$ with $|x_0 - y_0| =: r > 0$ and $B(x_0, 2\alpha r) \cup B(y_0, 2\alpha r) \subset \Omega$.*

Let $\phi : \mathbb{R}^{n \times n} \to \mathbb{R}$ be locally Lipschitz and such that

$$\phi(F) \geq \bar{c} \min \left\{ \mathrm{dist}^p(F, SO(n)U_1), |F| + 1 \right\}. \tag{2.1}$$

For every $u \in C^1(\Omega; \mathbb{R}^n)$ which satisfies

$$\frac{1}{r^n} \int_\Omega \phi(\nabla u) dx + \frac{1}{r^{n-1}} \int_\Omega |D[\phi(\nabla u)]| \leq \eta \tag{2.2}$$

there is a $\mathcal{L}^{2n}$-measurable set $E \subset B(x_0, \alpha r) \times B(y_0, \alpha r) =: B_1 \times B_2$ with $\mathcal{L}^{2n}(E) \geq (1-\theta)\mathcal{L}^{2n}(B_1 \times B_2)$ such that for all $(x, y) \in E$ one has

$$1 - c\varepsilon \leq \frac{|u(x) - u(y)|}{|U_1(x - y)|} \leq 1 + c\varepsilon, \tag{2.3}$$

where $[x, y]$ denotes the segment with endpoints $x$ and $y$, and

$$\varepsilon := \frac{1}{r^n} \int_\Omega \mathrm{dist}(\nabla u, K) \, dx. \tag{2.4}$$

*Proof.* We use conv to denote the convex hull and $c$ a generic constant that can depend on $K$, $n$, $\alpha$, $\theta$, $\bar{c}$, $p$. Without loss of generality we can assume $SO(n)U_i \cap SO(n)U_j = \emptyset$, $r = 1$ (by scaling) and $\Omega = \mathrm{conv}(B(x_0, 2\alpha) \cup B(y_0, 2\alpha))$ (by restricting $u$), $\eta \leq 1$.

*Step 1.* We show that

$$\int_\Omega \mathrm{dist}(\nabla u, SO(n)U_1) \, dx \leq c\eta^{1/p}. \tag{2.5}$$

Indeed, by (2.1) we have

$$\mathrm{dist}(F, SO(n)U_1) \leq c\phi(F) + c\phi^{1/p}(F).$$

After integration, using the embedding of $L^p$ into $L^1$ and that $\eta \leq 1$, we obtain

$$\begin{aligned}
\int_\Omega \mathrm{dist}(\nabla u, SO(n)U_1) \, dx &\leq c\|\phi(\nabla u)\|_{L^1} + c\|\phi^{1/p}(\nabla u)\|_{L^1} \leq c\eta + c\|\phi^{1/p}(\nabla u)\|_{L^p} \\
&\leq c\eta + c\|\phi(\nabla u)\|_{L^1}^{1/p} \leq c\eta^{1/p}.
\end{aligned}$$

*Step 2.* We construct a Lipschitz function $v$ which agrees with $u$ on a large set, using the truncation argument discussed in Proposition 2.2 below.

By assumption $\phi(F) = 0$ implies $F \in SO(n)U_1$ and in particular $\det F = \det U_1$. By continuity of the lower bound in (2.1) we can choose $\phi_0 > 0$ so that $\phi(F) \leq \phi_0$ implies

$$\mathrm{dist}(F, K) = \mathrm{dist}(F, SO(n)U_1) \quad \text{and} \quad \det F \geq \frac{1}{2} \det U_1. \tag{2.6}$$

Let $\Omega' := \mathrm{conv}\bigl(B(x_0, 15\alpha/8) \cup B(y_0, 15\alpha/8)\bigr)$. By Proposition 2.2, applied to $u$ on the pair of sets $\Omega' \subset\subset \Omega$ with this value of $\phi_0$, we obtain, for sufficiently small $\eta$, a set $\omega$ such that $\phi(\nabla u) \leq \phi_0$ on $\Omega' \setminus \omega$, and

$$\mathcal{L}^n(\omega) + \mathcal{H}^{n-1}(\partial\omega) \leq c\eta. \tag{2.7}$$

The set $\omega$ is the union of finitely many closed balls. Let $v : \mathbb{R}^n \to \mathbb{R}^n$ be a Lipschitz map which coincides with $u$ on $\Omega' \setminus \omega$. This exists by Kirszbraun's theorem, or by the – for the present purposes also sufficient – simpler result in [15, Sec. 3.1, Th. 1]. We conclude that

$$\mathrm{Lip}(v) \leq L \quad \text{and} \quad \text{for every } x \in \Omega' \setminus \omega \text{ the matrix } F = \nabla v(x) = \nabla u(x) \text{ obeys (2.6)}. \tag{2.8}$$

Here the Lipschitz constant $L$ on $v$ does not depend on $u$, $\eta$ and $\varepsilon$.

*Step 3.* We show that $v$ has degree 1 on a domain $D$ slightly smaller than $\Omega'$, and use this to prove invertibility.

Using first that $\mathrm{Lip}(v) \leq L$ and then (2.5) and (2.7) we have

$$\int_{\Omega'} \mathrm{dist}^n(\nabla v, SO(n)U_1)\, dx \leq (L+|U_1|)^{n-1} \int_{\Omega'\setminus\omega} \mathrm{dist}(\nabla u, SO(n)U_1)\, dx + \int_\omega (L+|U_1|)^n\, dx$$
$$\leq c\eta^{1/p} + c\mathcal{L}^n(\omega) \leq c\eta^{1/p}.$$

Thus, from the quantitative one-well rigidity estimate by Friesecke, James and Müller [18], it follows that there is $Q \in SO(n)$ such that

$$\int_{\Omega'} |\nabla v - QU_1|^n\, dx \leq c\eta^{1/p}.$$

By the Poincaré inequality there is $b \in \mathbb{R}^n$ such that $v$ is close to an affine map $A(x) = QU_1 x + b$, i.e.,

$$\int_{\Omega'} |\nabla v - QU_1|^n + |v - A|^n dx \leq c\eta^{1/p}. \tag{2.9}$$

Therefore by a choice argument based on Fubini there is $\alpha_1 \in (7\alpha/4, 15\alpha/8)$ such that, setting $D := \mathrm{conv}\big(B(x_0, \alpha_1) \cup B(y_0, \alpha_1)\big)$,

$$\int_{\partial D} |\nabla v - QU_1|^n + |v - A|^n d\mathcal{H}^{n_1} \leq c\eta^{1/p}.$$

By the embedding of $W^{1,n}(\partial D; \mathbb{R}^n)$ in $L^\infty(\partial D; \mathbb{R}^n)$ we have

$$|v(x) - A(x)| \leq c\eta^{1/p} \qquad \text{for all } x \in \partial D. \tag{2.10}$$

We shall now prove invertibility of $v$ with an argument based on the Brouwer degree, as in [9, Proposition 2.2]. We first assert that

$$\deg(v, D, z) = 1 \text{ for all } z \in A(D'), \tag{2.11}$$

where $D' := \mathrm{conv}\big(B(x_0, 3\alpha/2) \cup B(y_0, 3\alpha/2)\big)$. Indeed, consider the $C^0$ homotopy defined by $v_t(x) := tv(x) + (1-t)A(x)$ for $t \in [0,1]$. By (2.10) the maps $v_0 = A$ and $v_1 = v$ are uniformly close on $\partial D$, and by convexity (2.10) holds also for all $v_t$. But $\mathrm{dist}(D', \partial D) \geq \alpha/4$ and therefore, $|A(x) - A(y)| \geq |U_1^{-1}|\alpha/4$ for all $x \in \partial D$, $y \in D'$. This implies that for sufficiently small $\eta$

$$v_t(\partial D) \cap A(D') = \emptyset \qquad \text{for any } t \in [0,1]. \tag{2.12}$$

Therefore the Brouwer degree $\deg(v_t, D, z)$ does not depend on $t$ for any $z \in A(D')$ (see [16, Theorem 2.3(2)]). The affine map $A$ has degree 1, hence (2.11) follows.

In order to make use of (2.11) to obtain injectivity and to estimate the size of counterimages of sets, we assert that

$$1 \leq \#(v^{-1}(z) \cap D) \leq 1 + 2\#(v^{-1}(z) \cap D \cap \omega) \qquad \text{for all } z \in A(D'). \tag{2.13}$$

Indeed, the first inequality in (2.13) is obvious from (2.11). The second follows from (2.11) and the fact that for every $x \in D \setminus \omega$ one has that $v$ is $C^1$ in a neighborhood of $x$, and $\det \nabla v(x) > 0$. As a special case, we remark that

$$\text{for any } z \in A(D') \setminus v(\omega) \text{ there is exactly one } x \in D \cap v^{-1}(z). \tag{2.14}$$

Now consider the set

$$\omega_0 := \{x \in D : v(x) \in A(D') \cap v(\omega)\} = v^{-1}\big(A(D') \cap v(\omega)\big) \cap D.$$

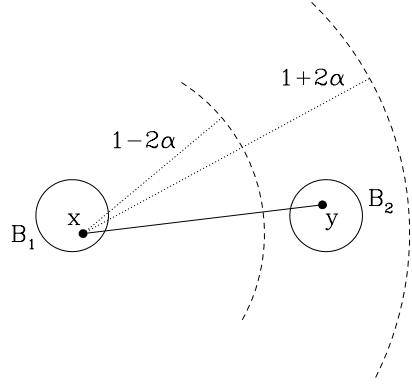

Figure 1: Integration domains in (2.17). For every given $x \in B_1$, the integration over $y$ is extended from $y \in B_2$ to the larger set $y \in B_{1+2\alpha}(x) \setminus B_{1-2\alpha}(x)$. Recall that $B_1$ and $B_2$ are balls with radius $\alpha$, the distance between the centers is 1.

To estimate its volume, we use $\mathcal{L}^n(\omega_0) \leq \mathcal{L}^n(\omega) + \mathcal{L}^n(\omega_0 \setminus \omega)$. In $\omega_0 \setminus \omega$ we know, by (2.6) and (2.8), that $\det \nabla v \geq \det U_1/2$. Therefore using the area formula we have

$$\frac{\det U_1}{2} \mathcal{L}^n(\omega_0 \setminus \omega) \leq \int_{\omega_0 \setminus \omega} \det \nabla v \, dx = \int_{\mathbb{R}^n} \#\bigl(v^{-1}(z) \cap \omega_0 \setminus \omega\bigr) dz.$$

By the definition of $\omega_0$ the latter integral can be restricted to $A(D') \cap v(\omega)$. Further, since $\omega_0 \subset D$, recalling (2.13) we obtain

$$\mathcal{L}^n(\omega_0 \setminus \omega) \leq c \int_{A(D') \cap v(\omega)} \#\bigl(v^{-1}(z) \cap D\bigr) dz \leq c \int_{A(D') \cap v(\omega)} 1 + \#\bigl(v^{-1}(z) \cap D \cap \omega\bigr) dz.$$

Finally, a second application of the area formula gives

$$\int_{v(\omega)} \#\bigl(v^{-1}(z) \cap \omega\bigr) dz = \int_\omega |\det \nabla v| \, dx \leq L^n \, \mathcal{L}^n(\omega),$$

and since $\mathcal{L}^n(v(\omega)) \leq L^n \mathcal{L}^n(\omega) \leq c\eta$ we conclude

$$\mathcal{L}^n(\omega_0) \leq \mathcal{L}^n(\omega) + \mathcal{L}^n(\omega_0 \setminus \omega) \leq c\eta. \tag{2.15}$$

*Step 4.* We finally turn to the construction of the "good" pairs $(x,y) \in B_1 \times B_2$. They have to be good in two aspects: first, the line integral of $\mathrm{dist}(\nabla u, SO(n)U_1)$ along the segment $[x,y]$ should be small, implying the upper bound in (2.3); second, there should be a curve $\gamma_{xy}: [0,1] \to D$ which is a parametrization of the preimage of $[v(x), v(y)]$ such that the integral of the same quantity along $\gamma_{xy}$ is also small, implying the lower bound in (2.3).

We say that a property $P$ holds for most choices of $(x,y)$ if the set where it does not hold has measure bounded by some $h(\eta)$, with $h(\eta) \to 0$ as $\eta \to 0$. Here $h$ should only depend on $K$, $n$, $\bar{c}$, $\alpha$, $\theta$, $p$, but not on $u$. For example, (2.15) shows that for most choices of $(x,y)$ one has $x, y \notin \omega_0$.

We start by proving that

$$\frac{1}{r} \int_{[x,y]} \mathrm{dist}(\nabla u, SO(n)U_1) \, d\mathcal{H}^1 \leq c\varepsilon \tag{2.16}$$

for most pairs $(x,y) \in B_1 \times B_2$. The key idea is to extend the integral to the entire line through $x$ and $y$, then integrate over $(x,y) \in B_1 \times B_2$, and change variables (see Figure 1). To simplify the tracking of the integration domains we set $f := \mathrm{dist}(\nabla u, K)$ in $\Omega$ and 0 elsewhere. We compute, with $z := y - x$,

$$\int_{B_1 \times B_2} \left( \int_{[x,y]} f \, d\mathcal{H}^1 \right) dx dy \leq \int_{B_1} \int_{1-2\alpha < |z| < 1+2\alpha} \int_0^1 |z| f(x+tz) \, dt \, dz \, dx. \qquad (2.17)$$

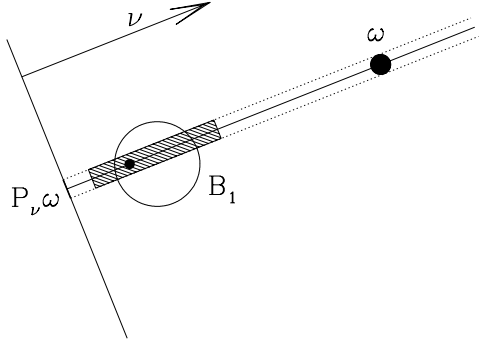

Figure 2: Geometry entering the definition of the set $G(\nu)$ in (2.18). Here $\omega$ is a single closed ball (black), its projection on $\nu^\perp$ an interval, and the dashed cylinder is the set on the right-hand side. The set $G_\nu$ is the intersection of the cylinder with $B_1$.

We swap the order of integration and exploit that for any $t \in \mathbb{R}$, $z \in \mathbb{R}^n$, one has

$$\int_{\mathbb{R}^n} f(x+tz) dx = \int_{\mathbb{R}^n} f(x) dx = \int_\Omega \mathrm{dist}(\nabla u, K) dx = \varepsilon.$$

Therefore

$$\int_{B_1 \times B_2} \left( \int_{[x,y]} f \, d\mathcal{H}^1 \right) dx dy \leq \int_{1-2\alpha < |z| < 1+2\alpha} |z| dz \int_0^1 dt \, \varepsilon \leq c\varepsilon.$$

Hence, $\int_{[x,y]} f \, d\mathcal{H}^1 \leq c\eta^{-1}\varepsilon$ for most pairs $(x,y)$, i.e., for all $(x,y) \in (B_1 \times B_2) \setminus I$, where $I$ is a certain set with $\mathcal{L}^{2n}(I) \leq \eta$.

Now we show that $[x,y] \cap \omega = \emptyset$ for most choices of $(x,y)$. By (2.6) and (2.8) this will imply that $f = \mathrm{dist}(\nabla u, SO(n)U_1)$ on $[x,y]$, and hence (2.16).

Fix a unit vector $\nu \in \mathbb{S}^{n-1}$, and let $P_\nu$ be the projection onto $\nu^\perp$. Since $\omega$ is a finite union of closed balls, $P_\nu \omega = P_\nu \partial\omega$. From $\mathcal{H}^{n-1}(P_\nu \omega) = \mathcal{H}^{n-1}(P_\nu \partial\omega) \leq \mathcal{H}^{n-1}(\partial\omega)$ we obtain that for each $\nu$ the set

$$G(\nu) := \{ x \in B_1 : x + \mathbb{R}\nu \cap \omega \neq \emptyset \} \subset P_\nu \omega \times [x_0 \cdot \nu - \alpha, x_0 \cdot \nu + \alpha]\nu \qquad (2.18)$$

has volume smaller then $c\eta$ (see Figure 2). Then

$$\int_{B_1 \times B_2} \chi_{[x,y] \cap \omega \neq \emptyset}(x,y) dy dx \leq \int_{B_1} \int_{1-2\alpha < |z| < 1+2\alpha} \chi_{x+\mathbb{R}z \cap \omega \neq \emptyset}(x,y) dz dx$$
$$\leq \int_{1-2\alpha < |z| < 1+2\alpha} \mathcal{L}^n(G(z/|z|)) \, dz \leq c\eta.$$

This concludes the proof of (2.16).

We now turn to the inverse, and assert that for most $(x,y) \in B_1 \times B_2$

$$[v(x), v(y)] \subset A(D') \quad \text{and} \quad [v(x), v(y)] \cap v(\partial \omega) = \emptyset. \tag{2.19}$$

To see this, we first observe that by (2.9) and $\text{dist}(B(x_0, \alpha), \partial D') > \alpha/2$ it follows that for most $(x, y)$ one has $v(x), v(y) \in A(D')$, and since $A(D')$ is convex this implies $[v(x), v(y)] \subset A(D')$. It remains to prove the second part of (2.19). Since $\mathcal{H}^{n-1}(v(\partial \omega)) \leq c\eta$, for any $\nu \in \mathbb{S}^{n-1}$ we have $\mathcal{H}^{n-1}(P_\nu v(\partial \omega)) \leq c\eta$, and arguing as above, since $v(B_2)$ is bounded, the set

$$H := \{(X, Y) \in v(B_1) \times v(B_2) : [X, Y] \cap v(\partial \omega) \neq \emptyset\}$$

is small, in the sense that $\mathcal{L}^{2n}(H) \leq c\eta$. Let now

$$h := \{(x, y) \in (B_1 \setminus \omega_0) \times (B_2 \setminus \omega_0) : [v(x), v(y)] \cap v(\partial \omega) \neq \emptyset\}.$$

We claim that $\mathcal{L}^{2n}(h) \leq c\eta$. Indeed, if $(x, y) \in h$ then $\det \nabla v(x)$ and $\det \nabla v(y)$ are larger than $\det U_1/2$, and using the area formula as above leads to

$$\mathcal{L}^{2n}(h) \leq c \int_H \#(v^{-1}(X) \cap B_1 \setminus \omega_0) \cdot \#(v^{-1}(Y) \cap B_2 \setminus \omega_0) d(X, Y) = c\mathcal{L}^{2n}(H) \leq c\eta$$

where we used that if $x \in B_1 \setminus \omega_0$ then $\#v^{-1}(v(x)) \cap D = 1$, and analogously for $y$. Recalling (2.7) and (2.15), the proof of (2.19) is concluded.

Pick now a pair $(x, y) \in (B_1 \setminus \omega \setminus \omega_0) \times (B_2 \setminus \omega \setminus \omega_0)$ such that (2.19) holds. We claim that there is a curve $\gamma_{xy} \in C^1([0, T]; D \setminus \omega)$ such that $\gamma(0) = x$, $\gamma(1) = y$, $u(\gamma_{xy}([0, T])) = [v(x), v(y)]$, $u \circ \gamma_{xy}$ injective. To prove this, we first observe that, by the continuity of $\nabla u$, $\det \nabla u \geq \det U_1/2$ also on $\overline{D \setminus \omega}$, and, by the implicit function theorem, for any point $a \in \overline{D \setminus \omega}$ there is $\rho_a > 0$ such that $u$ has a $C^1$ inverse as a map from $B(a, \rho_a)$ to $u(B(a, \rho_a))$. Notice that if $a \notin \partial D \cup \partial \omega$ then we can assume $B(a, \rho_a) \subset D \setminus \omega$ since $D \setminus \omega$ is open.

Let $\Gamma := v^{-1}([v(x), v(y)]) \cap \overline{D \setminus \omega}$. Since $\Gamma$ is compact, it can be covered by a finite number of balls $B(a, \rho_a)$ as above, $a \in \Gamma$. But, from $[v(x), v(y)] \subset A(D')$ and (2.12), we have $[v(x), v(y)] \cap v(\partial D) = \emptyset$; thus, recalling also (2.19), no such ball is centered on $\partial D \cup \partial \omega$. Therefore all those balls are contained in $D \setminus \omega$, and there $u = v$. This proves that $\Gamma$ is a finite union of $C^1$ arcs, which do not touch $\partial D \cup \partial \omega$. Therefore they can have endpoints only in the sets $v^{-1}(v(x))$ and $v^{-1}(v(y))$. But since $x, y \notin \omega_0$, the points $v(x)$ and $v(y)$ have only one counterimage each, namely, $x$ and $y$. Hence $\Gamma$ consists of a single arc, with endpoints $x$ and $y$.

We show that for most $(x, y)$ the curve $\gamma_{xy}$ carries energy of order $\varepsilon$. Define $g : \mathbb{R}^n \to \mathbb{R}$ as

$$g(z) := \sum_{x \in u^{-1}(z) \cap (D \setminus \omega)} \text{dist}(\nabla u(x), K) = \sum_{x \in u^{-1}(z) \cap (D \setminus \omega)} \text{dist}(\nabla u(x), SO(n)U_1),$$

where it is understood that $g = 0$ if the sum is empty, and the equality follows from (2.6), (2.8). By a change of variables we have

$$\int_{\mathbb{R}^n} g(z) dz = \int_{D \setminus \omega} f(x) \det \nabla u(x) dx \leq c\varepsilon.$$

Arguing as in (2.17) and the following, we see that the set $M$ of $(\xi, \zeta) \in A(B(x_0, 3\alpha/2)) \times A(B(y_0, 3\alpha/2))$ where $\int_{[\xi, \zeta]} g \, d\mathcal{H}^1 \leq c\varepsilon/\eta$ does not hold is small. Then, so is the set $m$ of $(x, y) \in (B_1 \setminus \omega_0) \times (B_2 \setminus \omega_0)$ where the condition is violated for $(\xi, \zeta) = (u(x), u(y))$ (this argument is just as for the case of the two sets $H$ and $h$ above). Hence for most pairs $(x, y)$

$$\int_{\gamma_{xy}} \text{dist}(\nabla u, SO(n)U_1) \, d\mathcal{H}^1 \leq c \int_{[u(x), u(y)]} g \, d\mathcal{H}^1 \leq c_\eta \varepsilon \,.. \tag{2.20}$$

*Step 5.* We show that, from Step 4, (2.3) follows. We work with a pair $(x,y)$ satisfying (2.16) and (2.20). Using (2.16), we have

$$\begin{aligned}|u(x)-u(y)| &\leq \int_{[x,y]}|\nabla_\tau u|d\mathcal{H}^1 \leq |U_1(x-y)| + \int_{[x,y]} \mathrm{dist}(\nabla u, SO(n)U_1)\,d\mathcal{H}^1 \\ &\leq |U_1(x-y)| + c\varepsilon\,,\end{aligned} \qquad (2.21)$$

where $\nabla_\tau u$ denotes the tangential gradient to $[x,y]$.

To prove the converse inequality, we use the fact that there is a curve $\gamma \in C^1([0,T]; D\setminus\omega)$ such that $|\gamma'|=1$, $\gamma(\{0,1\}) = \{x,y\}$, and $u\circ\gamma$ is monotone parametrization of the segment $[u(x), u(y)]$. This implies in particular that the directional derivative $\nabla u(\gamma)\gamma'$ is parallel to $u(x)-u(y)$, and hence that

$$|u(x)-u(y)| = \left|\int_0^T \frac{d}{dt} u(\gamma(t))dt\right| = \int_0^T |\nabla u(\gamma)\gamma'|\,dt\,.$$

For each $t\in [0,T]$, let $Q(t)$ be such that $\mathrm{dist}(\nabla u(\gamma(t)), SO(n)U_1) = |\nabla u(\gamma(t)) - Q(t)U_1|$. Then we obtain

$$\begin{aligned}\int_0^T |\nabla u(\gamma)\gamma'|\,dt &\geq \int_0^T |Q(t)U_1\gamma'(t)| - |(\nabla u(\gamma(t)) - Q(t)U_1)\gamma'(t)|\,dt \\ &\geq \int_0^T |U_1\gamma'|\,dt - \int_0^T |\nabla u\circ\gamma - QU_1|\,dt \\ &\geq \left|\int_0^T U_1\gamma'\,dt\right| - \int_0^T \mathrm{dist}(\nabla u\circ\gamma, SO(n)U_1)dt \\ &= |U_1(x-y)| - \int_{\gamma([0,T])} \mathrm{dist}(\nabla u, SO(n)U_1)d\mathcal{H}^1 \geq |U_1(x-y)| - c\varepsilon\,.\end{aligned}$$

This concludes the proof. □

We finally present the truncation that was used in Step 2. Precisely, we show that $u$ is $c$-Lipschitz on a subdomain $\Omega'$ compactly contained in $\Omega$, away from a set $\omega \subset \Omega$ with small volume and perimeter. Then it is possible to obtain a Lipschitz extension to $\mathbb{R}^n$, which coincides with the old function on $\Omega'\setminus\omega$. With respect to the classical construction discussed in [15, Sect. 6.6.2], we need here to estimate also the perimeter of $\omega$, hence need to construct, via a covering argument, a smoother set (a countable union of balls).

**Proposition 2.2.** *Let $\Omega'\subset\subset\Omega\subset\mathbb{R}^n$, both open, with $\Omega$ bounded and $\Omega'$ connected and Lipschitz. Let $\phi\in\mathrm{Lip}_{\mathrm{loc}}(\mathbb{R}^{n\times n};[0,\infty))$ be such that*

$$\phi(F) \geq \frac{1}{C_0}|F| - C_0 \qquad (2.22)$$

*for some $C_0 > 0$, and $\phi_0$ be a positive constant. Then there is $M > 0$ (depending on $C_0$, $\phi_0$, $\Omega$ and $\Omega'$) such that, for any $u\in W^{1,1}(\Omega;\mathbb{R}^n)$ there is a set $\omega$ with*

$$\mathcal{L}^n(\omega) + \mathrm{Per}_\Omega(\omega) \leq M\left[\int_\Omega \phi\circ\nabla u\,dx + \int_\Omega |D(\phi\circ\nabla u)|\right] \qquad (2.23)$$

*such that the map $u$ is $M$-Lipschitz on $\Omega'\setminus\omega$ and $\phi(\nabla u) \leq \phi_0$ on $\Omega'\setminus\omega$.*

*If additionally $u\in C^1(\Omega;\mathbb{R}^n)$, then $\omega$ can be chosen to be the union of finitely many closed balls.*

*Proof.* We set

$$\eta := \int_\Omega \phi\circ\nabla u\,dx + \int_\Omega |D(\phi\circ\nabla u)| \leq \eta_0\,. \qquad (2.24)$$

We first show that there is $c_0 \in (\phi_0/3, 2\phi_0/3)$ such that the set

$$\omega' := \{x \in \Omega : \phi(\nabla u(x)) \geq c_0\} \qquad (2.25)$$

obeys $\mathcal{L}^n(\omega') + \mathrm{Per}_\Omega(\omega') \leq c\eta$. Indeed, by the coarea formula for BV functions [1, Th. 3.40] and from (2.24) we have $\int_{\mathbb{R}} \mathrm{Per}_\Omega(\Omega_t)\, dt = |D[\phi(\nabla u)]|(\Omega) \leq \eta$, where $\Omega_t := \{x \in \Omega : \phi(\nabla u(x)) \geq t\}$. Therefore there is $c_0 \in (\phi_0/3, 2\phi_0/3)$ such that the set $\omega' = \Omega_{c_0}$ satisfies $\mathrm{Per}_\Omega(\omega') \leq 3\eta/\phi_0$. By (2.24) we get also $\mathcal{L}^n(\omega') \leq 3\eta/\phi_0$.

Let $r_0 := \mathrm{dist}(\Omega', \partial\Omega)$, for $h \in \mathbb{N}$ set $r_h := r_0 2^{-h}$, and for each $x \in \Omega'$ consider the sequence of balls $B_h(x) := B(x, r_h)$. We say that a ball $B = B_h(x)$ is good if

$$\int_B |D(\phi \circ \nabla u)| \leq \mathcal{H}^{n-1}(\partial B) \quad \text{and} \quad \mathcal{L}^n(B \cap \omega') \leq \frac{1}{2^{n+1}} \mathcal{L}^n(B). \qquad (2.26)$$

By assumption, all the level-0 balls $B_0(x)$ are compactly contained in $\Omega$. If there is one $x$ such that $B_0(x)$ is not good, then $\eta \geq c$, and the conclusion follows with $\omega$ a large ball containing $\Omega$. Therefore we can assume that all balls $B_0(x)$, $x \in \Omega'$, are good. Let $\omega'' \subset \Omega'$ be the set of all centers of bad balls. For each $x \in \omega''$ let $h(x)$ be the smallest $h$ such that $B_h(x)$ is bad. We have already shown that $h(x) > 0$ for all $x \in \omega''$. Let $\mathcal{F}$ be the family of closed balls $\overline{B}_{h(x)}(x)$, for $x \in \omega''$. By the Besicovitch covering theorem, $\omega''$ can be covered by a fixed number $N_B$ of disjoint families of such balls, let $A_1...A_{N_B}$ be the sets of their centers. Let

$$\omega_k := \bigcup_{x \in A_k} \overline{B}_{h(x)}(x) \qquad (2.27)$$

and

$$\omega := \{x \in \Omega : x \text{ is not a Lebesgue point for } u \text{ or } \nabla u\} \cup \bigcup_{k=1}^{N_B} \omega_k.$$

Obviously $\omega'' \subset \cup_k \omega_k \subset \omega$. We now show that $\phi(\nabla u) \leq c_0$ on $\Omega' \setminus \omega$. Indeed, if we had $x \in \Omega' \setminus \omega \subset \Omega' \setminus \omega''$ and $\phi(\nabla u(x)) > c_0$ then by continuity of $\phi$ there would be $\varepsilon > 0$ such that $\phi(F) > c_0$ for all $F \in \mathbb{R}^{n \times n}$ with $|F - \nabla u(x)| < \varepsilon$. In particular, this would imply that $|\nabla u(y) - \nabla u(x)| \geq \varepsilon$ for all $y$ such that $\phi(\nabla u(y)) < c_0$, i.e., for all $y \in \Omega \setminus \omega'$. Then

$$\frac{1}{\mathcal{L}^n(B(x, r_h))} \int_{B(x, r_h)} |\nabla u(y) - \nabla u(x)|\, dy \geq \varepsilon \frac{\mathcal{L}^n(B(x, r_h) \setminus \omega')}{\mathcal{L}^n(B(x, r_h))}.$$

Since $x$ is a Lebesgue point for $\nabla u$, the left-hand side converges to zero as $h \to \infty$, which implies that

$$\lim_{h \to \infty} \frac{\mathcal{L}^n(B(x, r_h) \setminus \omega')}{\mathcal{L}^n(B(x, r_h))} = 0,$$

Therefore for sufficiently large $h$ the ball $B(x, r_h)$ cannot be good, that is, $x \in \omega''$, which is a contradiction. We conclude that $\phi(\nabla u) \leq c_0$ on $\Omega' \setminus \omega$.

We shall prove that each $\omega_k$ obeys the volume and perimeter estimate, the same will then hold for their union (with a different constant).

We consider $\omega_k$, which is the union of countably many disjoint bad balls, such that for each of them the ball twice as large is good (by the minimality of $h(x)$, and the fact that $h(x) \neq 0$). For each of those balls, we claim that

$$\mathcal{H}^{n-1}(\partial B) \leq \int_B |D(\phi \circ \nabla u)| + c \, \mathrm{Per}_B(\omega'). \qquad (2.28)$$

If (2.28) holds, by $\mathcal{L}^n(B) \leq r_0 \mathcal{H}^{n-1}(\partial B)$ an analogous estimate holds for the volume. Further, since the balls composing $\omega_k$ are disjoint, we obtain

$$\begin{aligned}
\mathcal{L}^n(\omega_k) + \operatorname{Per}_\Omega(\omega_k) &= \sum_{x \in A_k} \mathcal{L}^n(B_{h(x)}(x)) + \mathcal{H}^{n-1}(\partial B_{h(x)}(x)) \\
&\leq \sum_{x \in A_k} c \int_{B_{h(x)}(x)} |D(\phi \circ \nabla u)| + c \operatorname{Per}_{B_{h(x)}(x)}(\omega') \\
&\leq c \int_\Omega |D(\phi \circ \nabla u)| + c \operatorname{Per}_\Omega(\omega') \leq c\eta \,.
\end{aligned}$$

Summing over $k$ from 1 to $N_B$ gives (2.23).

We now prove (2.28). Recall that $B$ is a bad ball contained in a good ball $B'$ twice as large. Since $B$ is bad, one of the two conditions in (2.26) is violated. If it is the first one, (2.28) follows immediately. If it is the second one, using the fact that $B'$ is good we obtain

$$\frac{1}{2^{n+1}} \mathcal{L}^n(B) \leq \mathcal{L}^n(B \cap \omega') \leq \mathcal{L}^n(B' \cap \omega') \leq \frac{1}{2^{n+1}} \mathcal{L}^n(B') = \frac{1}{2} \mathcal{L}^n(B) \,.$$

Therefore both $\mathcal{L}^n(B \cap \omega')$ and $\mathcal{L}^n(B \setminus \omega')$ have volume larger than $2^{-(n+1)} \mathcal{L}^n(B)$, and by the relative isoperimetric inequality [1, Eq. (3.37)] we obtain $\operatorname{Per}_B(\omega') > c \mathcal{H}^{n-1}(\partial B)$. This concludes the proof of (2.28).

It remains to show that $u$ is Lipschitz on $\Omega' \setminus \omega$. Pick $x \in \Omega' \setminus \omega$. Then $x \notin \omega''$, hence for every $h \in \mathbb{N}$ the ball $B_h(x)$ is good. By the first condition in (2.26) and the Poincaré inequality there is $\widehat{\phi} = \widehat{\phi}(x,h) \in \mathbb{R}$ such that

$$\int_{B_h(x)} |\phi(\nabla u) - \widehat{\phi}| dy \leq c \mathcal{L}^n(B_h) \,, \tag{2.29}$$

and by the second one $|\widehat{\phi}| \leq c$. Therefore, using (2.22),

$$\int_{B_h(x)} |\nabla u| dy \leq C_0 \int_{B_h(x)} (\phi(\nabla u) + C_0) \, dy \leq c \mathcal{L}^n(B_h) \,, \tag{2.30}$$

for every $h$. Let

$$F_h(x) := \frac{1}{\mathcal{L}^n(B_h)} \int_{B_h(x)} u(y) dy \,.$$

By the Poincaré inequality, (2.30), and the fact that $B_{h+1}(x) \subset B_h(x)$, we obtain

$$\begin{aligned}
|F_{h+1}(x) - F_h(x)| &= \left| \frac{1}{\mathcal{L}^n(B_{h+1})} \int_{B_{h+1}(x)} (u(y) - F_h(x)) dy \right| \\
&\leq \frac{1}{\mathcal{L}^n(B_{h+1})} \int_{B_{h+1}(x)} |u(y) - F_h(x)| dy \\
&\leq \frac{1}{\mathcal{L}^n(B_{h+1})} \int_{B_h(x)} |u(y) - F_h(x)| dy = \frac{2^n}{\mathcal{L}^n(B_h)} \int_{B_h(x)} |u(y) - F_h(x)| dy \\
&\leq c r_h \frac{1}{\mathcal{L}^n(B_h)} \int_{B_h(x)} |\nabla u| \leq c r_h \,,
\end{aligned}$$

with a universal constant $c$. Therefore, summing the geometric series and using that $\lim_{h \to \infty} F_h(x) = u(x)$, we obtain

$$|u(x) - F_h(x)| \leq c r_h \tag{2.31}$$

for all $h$. Arguing analogoulsy with the pair $B_0(x) \subset \Omega$ we also obtain

$$|u(x) - F| \leq c \tag{2.32}$$

where $F := \mathcal{L}^n(\Omega)^{-1} \int_\Omega u(y) dy$.

We now show that for all $x, y \in \Omega' \setminus \omega$ one has

$$|u(x) - u(y)| \leq c'|x - y|. \tag{2.33}$$

To prove (2.33), assume first that $r_0 \geq 2|x-y|$, and let $k$ be the largest integer such that $r_k \geq 2|x-y|$. By maximality of $k$, $r_k \leq 4|x-y|$. This implies that the balls $B_k(x)$ and $B_k(y)$ have a large overlap, and arguing with Poincaré's inequality as above, from (2.30) we obtain

$$|F_k(x) - F_k(y)| \leq c r_k.$$

Therefore, recalling (2.31),

$$|u(x) - u(y)| \leq |u(x) - F_k(x)| + |F_k(x) - F_k(y)| + |F_k(y) - u(y)| \leq 3cr_k \leq c'|x-y|$$

which concludes the proof of (2.33) in the case $r_0 \geq 2|x-y|$.

We now turn to the case $2|x-y| \geq r_0$. Then arguing analogously on the basis of (2.32) gives

$$|u(x) - u(y)| \leq |u(x) - F| + |F - u(y)| \leq 2c \leq c'|x-y|$$

This concludes the proof of (2.33) and hence of the first statement of the proposition.

Assume now that $\nabla u$, and hence $\phi(\nabla u)$, are continuous. Possibly restricting to a smaller set we can assume them to be uniformly continuous. We briefly discuss the few changes necessary to the construction in this case. Let $j > 0$ be such that $|x - y| \leq 2r_j = 2r_0 2^{-j}$ implies $|\phi(\nabla u(x)) - \phi(\nabla u(y))| \leq \phi_0/3$. We say that a ball $B(x, r_h)$ is good if either (2.26) holds or $h > j$ and $B(x, r_j)$ is good. Let $\omega'' \subset \Omega'$ be again the set of centers of bad balls, and for each $x \in \omega''$ we let $h(x)$ be the smallest $h \in \mathbb{N}$ such that $B_h(x)$ is bad. Then $h(x) \leq j$ for all $x \in \omega''$, hence proceeding with Besicovitch as above, each $\omega_k \subset \Omega$ is composed by only finitely many balls. Analogously, all points of $\Omega$ are Lebesgue points for the continuous functions $u$ and $\nabla u$, hence $\omega$ is a finite union of closed balls.

We now prove that this exceptional set also suffices. The proof of (2.28) is unchanged. Analogously, (2.30) hold for all $h \leq j$. In order to prove the same estimate for $h > j$, fix some $x \notin \omega''$. Since $B_j(x)$ is good there is $y \in B_j(x) \setminus \omega'$, which by definition of $\omega'$ obeys $\phi(\nabla u(y)) \leq 2\phi_0/3$. For any $z \in B_j(x)$ one has $|z - y| \leq \text{diam}(B_j) = 2r_j$, and therefore

$$\text{for all } z \in B_j(x) \text{ one has } \phi(\nabla u(z)) \leq \phi(\nabla u(y)) + \frac{1}{3}\phi_0 \leq \phi_0. \tag{2.34}$$

This uniform estimate proves (2.30) for all $h \geq j$. Further, it proves that $\phi(\nabla u) \leq \phi_0$ on $\Omega \setminus \omega \subset \Omega \setminus \omega''$. The proof that (2.30) implies that $u$ is $M$-Lipschitz on $\Omega' \setminus \omega$ is unchanged. □

## 3 Multiwell rigidity in $L^1$

In this section we use the segment rigidity of the last section to prove Theorem 1.1. The key step is a geometric argument permitting to prove the statement for the case that $D$ and $D'$ are two concentric balls, with $D'$ much smaller than $D$ (Proposition 3.4, based on Lemma 3.5, which in turn is based on Proposition 2.1, Lemma 3.2 and Lemma 3.3). Then a covering argument leads to Theorem 1.1. We start by three lemmas dealing with $n$-dimensional simplexes. We denote by $(e_h)_h$ the canonical basis of $\mathbb{R}^n$.

**Lemma 3.1.** *Let $\{a_1, a_2 \ldots a_{n+1}\} \subset \mathbb{R}^n$, $n \geq 2$, be such that the $n$ vectors $a_h - a_{n+1}$, for $h = 1 \ldots n$, are linearly independent. Let $\bar{a} := \sum_{h=1}^{n+1} a_h/(n+1)$ be the barycenter. Then there are $\rho$ and $C > 0$ (depending on the $a_h$ and on $n$) such that the following holds:*

*For all choices $\{b_1, b_2 \ldots b_{n+1}\} \subset \mathbb{R}^n$ with $b_h \in B(a_h, \rho)$ and all $p \in B(\bar{a}, \rho)$ there are $c_1, c_2, \ldots c_{n+1} \in (1/C, C)$ such that, for all $q \in \mathbb{R}^n$, one has*

$$\sum_{h=1}^{n+1} c_h |p - b_h| \leq \sum_{h=1}^{n+1} c_h |q - b_h|. \tag{3.1}$$

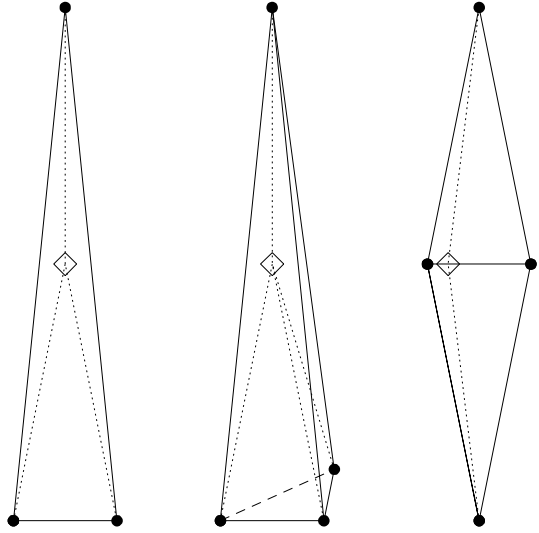

Figure 3: Left panel: sketch of the geometry in the proof of Proposition 3.4 for the two-dimensional case. The edges of the simplex (triangle) are rigid in the sense of Proposition 2.1, hence the map $u$ is close to the affine map $A$ on the vertices. The result is obtained estimating the length of the segments joining the vertices with a generic point $P$ close to the barycenter. Central panel: corresponding sketch in three dimensions. Right panel: for a comparison, we report the geometry used in [9]. There are five rigid connections among four points, and $P$ belongs to the central segment, which is rigid by (2.16). A direct extension of this geometry to higher dimension fails because the "mid-segment" becomes a "mid-plane", which is not any more rigid.

*Proof.* Since the vectors $\{a_h - a_{n+1}\}$ are linearly independent, the set

$$T_a := \mathrm{conv}\{a_1 \ldots a_{n+1}\} = \left\{\sum_{h=1}^{n+1} \lambda_h a_h : \lambda_h \in [0,1], \sum_{h=1}^{n+1} \lambda_h = 1\right\}$$

is a non-degenerate simplex in $\mathbb{R}^n$, with positive $n$-dimensional volume $\mathcal{L}^n(T_a) = |\det(\sum_{h=1}^n (a_h - a_{n+1}) \otimes e_h)|$. In particular the point $\bar{a}$ belongs to the interior of $T$. Let $\rho > 0$ be such that

$$B(\bar{a}, 3\rho) \subset T_a \,.$$

This implies that, for all admissible choice of the $b_h$'s,

$$B(\bar{a}, 2\rho) \subset T_b := \mathrm{conv}\{b_1 \ldots b_{n+1}\} \,. \tag{3.2}$$

To see this, notice that any point $x \in \partial T_b$ can be written as $x = \sum \mu_h b_h$, with $\mu_h \in [0,1]$, $\sum \mu_h = 1$, and at least one vanishing. Consider the point $y := \sum \mu_h a_h$. Clearly $y \in \partial T_a$, hence $|y - \bar{a}| \geq 3\rho$. At the same time $|x - y| \leq \sum \mu_h |a_h - b_h| < \rho \sum \mu_h = \rho$. Therefore $|x - \bar{a}| \geq |y - \bar{a}| - |x - y| > 2\rho$, and (3.2) follows.

Consider now any $p \in B(\bar{a}, \rho) \subset T_b$. Let $\lambda_h \in [0,1]$ with $\sum_{h=1}^{n+1} \lambda_h = 1$ be such that

$$p = \sum_{h=1}^{n+1} \lambda_h b_h \,.$$

We show that $\lambda_h \geq \rho/c$, where $c$ is a constant which depends only on the $a_h$. Consider for definiteness $\lambda_1$ (the others are treated the same way), assume $\lambda_1 \leq 1/2$, and set $\lambda_h^* := \lambda_h/(1-\lambda_1)$ and $p^* := \sum_{h \neq 1} \lambda_h^* b_h \in \partial T_b$. From $|\lambda_h - \lambda_h^*| = \lambda_1 \lambda_h/(1-\lambda_1) \leq 2\lambda_1 \lambda_h$ we obtain

$$|p - p^*| = \left|\lambda_1 b_1 + \sum_{h \neq 1}(\lambda_h - \lambda_h^*)b_h\right| \leq 3 \max_h |b_h| \lambda_1 \leq c\lambda_1 \,,$$

where $c = 3\max_h |a_h| + 3\rho$. But since $p^* \in \partial T_b$, (3.2) implies that $p^* \notin B(\bar{a}, 2\rho)$. Therefore
$$2\rho < |p^* - \bar{a}| < |p^* - p| + |p - \bar{a}| \leq c\lambda_1 + \rho\,,$$
which implies $\lambda_1 > \rho/c$. Analogously, $b_h \in \partial T_b$ implies $b_h \notin B(\bar{a}, 2\rho)$, and therefore $|p - b_h| \geq \rho$ for all $h$.

We set $c_h := \lambda_h |p - b_h|$. The above argument shows that $1/C \leq c_h \leq C$ for some constant $C > 0$ depending only on the $a_h$. We define the function $f_p : \mathbb{R}^n \to \mathbb{R}$ by
$$f_p(q) := \sum_{h=1}^{n+1} c_h |q - b_h|\,.$$
We compute (for $q \notin \{b_h\}_h$)
$$\nabla f_p(q) = \sum_{h=1}^{n+1} c_h \frac{q - b_h}{|q - b_h|}\,,$$
and observe that
$$\nabla f_p(p) = \sum_{h=1}^{n+1} c_h \frac{p - b_h}{|p - b_h|} = \sum_{h=1}^{n+1} \lambda_h (p - b_h) = p - \sum_{h=1}^{n+1} \lambda_h b_h = 0\,.$$
Since $f_p$ is convex, this implies that $p$ is the global minimum of $f_p$, i.e., that
$$f_p(p) \leq f_p(q) \qquad \text{for all } q \in \mathbb{R}^n\,,$$
which concludes the proof. $\square$

**Lemma 3.2.** *Let $\delta > 0$, $n \in \mathbb{N}$, $n \geq 2$, be given. Then there are $\rho > 0$, $\mu > 0$, $c > 0$ depending only on $\delta$ and $n$ such that the following holds. Let $\{a_1, a_2 \ldots a_{n+1}\} = \{0, \mu e_1, \mu e_2, \ldots \mu e_{n-1}, e_n\} \subset \mathbb{R}^n$, and let $\bar{a} := \sum_{h=1}^{n+1} a_h/(n+1)$ be the barycenter.*
*For all choices $\{b_1, b_2 \ldots b_{n+1}\} \subset \mathbb{R}^n$ with $b_h \in B(a_h, \rho)$, and all $p \in B(\bar{a}, \rho)$, one has*
$$\left| \frac{p - b_h}{|p - b_h|} - e_n \right| \leq \delta\,. \tag{3.3}$$
*Further, there are $c_1, c_2, \ldots c_{n+1} \in (1/c, c)$ (depending on the $b_h$) such that for all $p \in B(\bar{a}, \rho)$, and all $q \in \mathbb{R}^n$, one has*
$$\sum_{h=1}^{n+1} c_h |p - b_h| \leq \sum_{h=1}^{n+1} c_h |q - b_h|\,. \tag{3.4}$$

*Proof.* To prove (3.3), we write for $h \leq n$
$$\left| p - b_h - \frac{1}{n+1} e_n \right| \leq |p - \bar{a}| + |a_h - b_h| + |a_h| + \left| \bar{a} - \frac{1}{n+1} e_n \right| \leq 2\rho + 2\mu\,,$$
and in the last case, since $a_{n+1} = e_n$,
$$\left| p - b_{n+1} + \frac{n}{n+1} e_n \right| \leq |p - \bar{a}| + |a_{n+1} - b_{n+1}| + \left| \bar{a} - e_n + \frac{n}{n+1} e_n \right| \leq 2\rho + \mu\,.$$
Therefore (3.3) holds provided $\mu$ and $\rho$ are chosen smaller than $c\delta$, with $c$ depending on the dimension $n$.

Fix one such $\mu$. Then the $a_h$'s are given, depending only on $\mu$, and satisfy the assumptions of Lemma 3.1. Hence there are $\rho$ and $c$, depending only on the $a_h$'s (and hence only on $\delta$) such that (3.4) holds. $\square$

The next lemma proves the well-known fact that if two simplexes have sides of approximately equal length, then one is close to an isometric copy of the other. For completeness we give a short self-contained proof.

**Lemma 3.3.** *Let $c > 0$, $a_1, \ldots a_{n+1} \in \mathbb{R}^n$ be such that $|a_h| \leq c$ for all $h$, and $\det(\sum_{h=1}^n (a_h - a_{n+1}) \otimes e_h) \geq 1/c$. Then for all $\varepsilon \in (0, 1/2)$ and all $b_1, \ldots b_{n+1} \in \mathbb{R}^n$ such that*

$$1 - \varepsilon \leq \frac{|b_h - b_k|}{|a_h - a_k|} \leq 1 + \varepsilon \tag{3.5}$$

*there are $R \in O(n) = \{R \in \mathbb{R}^{n \times n} : R^T R = \mathrm{Id}\}$ and $d \in \mathbb{R}^n$ such that for all $h$*

$$|b_h - (Ra_h + d)| \leq C\varepsilon$$

*where $C$ depends only on $c$ and $n$.*

*Proof.* By replacing the vectors with $\{a_h - a_{n+1}\}_h$ and $\{b_h - b_{n+1}\}_h$ we can assume $a_{n+1} = b_{n+1} = d = 0$. Let $F \in \mathbb{R}^{n \times n}$ be the unique linear map which satisfies $Fa_h = b_h$, for $h = 1, \ldots, n$ (this exists by the assumption on the determinant). Let $G := \sqrt{F^T F}$, and $R \in O(n)$ be such that $F = RG$. We need to show that $|G - \mathrm{Id}| \leq C\varepsilon$.

From (3.5), with $k = n + 1$, we obtain

$$|b_h^2 - a_h^2| \leq C\varepsilon$$

and

$$|(b_h - b_k)^2 - (a_h - a_k)^2| = |(b_h^2 - a_h^2) + (b_k^2 - a_k^2) - 2(b_h \cdot b_k - a_h \cdot a_k)| \leq C\varepsilon,$$

which imply

$$|b_h \cdot b_k - a_h \cdot a_k| \leq C\varepsilon,$$

where the constant $C$ may change from line to line and $a^2$ denotes the squared norm of the vector $a$. Consider two vectors $v, w \in \mathbb{S}^n$, and let $\lambda, \mu \in \mathbb{R}^n$ be such that $v = \sum_{h=1}^n \lambda_h a_h$, $w = \sum_{h=1}^n \mu_h a_h$. By the assumption on the determinant, $|\lambda| + |\mu| \leq C$. We compute

$$|v(F^T F - \mathrm{Id})w| = |(Fv) \cdot Fw - v \cdot w| = \left| \sum_{h,k} \lambda_h \mu_k (b_h \cdot b_k - a_h \cdot a_k) \right| \leq C\varepsilon \sum_{h,k} |\lambda_h \mu_k| \leq C\varepsilon.$$

This proves that $|F^T F - \mathrm{Id}| \leq C\varepsilon$. Therefore all eigenvalues of the positive-definite, symmetric matrix $F^T F$ lie in $[1 - C\varepsilon, 1 + C\varepsilon]$, and the same (with a different $C$) holds for $G$. This concludes the proof. □

This concludes the preparatory part, and we now come to the multiwell rigidity result in $L^1$. The rest of this section is organized as follows: we now state the result for balls (Prop. 3.4), then show how one can reduce to a special case, which is then proven in Lemma 3.5. Then we discuss how Proposition 3.4 implies Theorem 1.1. Finally, Remarks 3.6, 3.7 and 3.8 concern the fact that the separation assumption (1.5) on the matrices is satisfied for most cases of physical interest.

**Proposition 3.4.** *Let $U_1 \ldots U_l$ be matrices in $\mathbb{R}^{n \times n}$ with positive determinant, which obey (1.5). Let $K$ be defined as in (1.3). Then there are positive numbers $\eta$, $c^*$ and $\rho \in (0, 1)$ depending only on $\{U_i\}_i$ and $n$ such that for any ball $B := B(x_0, r) \subset \mathbb{R}^n$, any $u \in W^{1,1}(B; \mathbb{R}^n)$ with $\nabla u \in BV$ satisfying*

$$\frac{1}{r^{n-1}} \int_{B(x_0, r)} |D^2 u| \leq \eta \tag{3.6}$$

*one has*

$$\min_{J \in \{U_i\}_i} \int_{B(x_0, \rho r)} \mathrm{dist}(\nabla u, SO(n)J) \, dx \leq c^* \int_{B(x_0, r)} \mathrm{dist}(\nabla u, K) \, dx.$$

We divide the proof of Proposition 3.4 in two parts: we first reduce to a special case where one phase, corresponding to $SO(n)$, is dominant, and then treat that special case in Lemma 3.5 below.

*Proof.* By scaling we can assume $r = 1$, by density we can assume that $u$ is smooth, and by relabeling the wells that $U_1$ is the majority phase, i.e., that the set
$$E_1 := \{x \in B : \text{dist}(\nabla u(x), K) = \text{dist}(\nabla u(x), SO(n)U_1)\}$$
obeys $\mathcal{L}^n(E_1) \geq \mathcal{L}^n(B)/l$. This implies
$$\int_{E_1} \text{dist}(\nabla u, SO(n)U_1)\, dx = \int_{E_1} \text{dist}(\nabla u, K)\, dx \leq \int_B \text{dist}(\nabla u, K)\, dx =: \varepsilon, \qquad (3.7)$$
where the last equality defines $\varepsilon$. By the Poincaré inequality, (3.6) implies that there is $F \in \mathbb{R}^{n \times n}$ such that
$$\int_B |\nabla u - F|\, dx \leq c\eta. \qquad (3.8)$$
Thus, by (3.7) and (3.8),
$$\mathcal{L}^n(E_1) \text{dist}(F, SO(n)U_1) \leq \int_{E_1} \text{dist}(\nabla u, SO(n)U_1)\, dx + \int_B |\nabla u - F|\, dx \leq \varepsilon + c\eta,$$
and therefore
$$\int_B \text{dist}(\nabla u, SO(n)U_1)\, dx \leq \int_B \text{dist}(F, SO(n)U_1)\, dx + \int_B |\nabla u - F|\, dx \leq c(\varepsilon + \eta). \qquad (3.9)$$
In the case $\varepsilon \geq \eta$ the proof is concluded. In the following we can therefore assume that $\varepsilon \leq \eta$.

We fix a vector $\xi_1$ as in (1.5), a matrix $Q \in SO(n)$ such that $U_1\xi_1 = Qe_n$, and define
$$\tilde{u}(x) := u(U_1^{-1}Qx), \qquad \tilde{U}_i := U_i U_1^{-1} Q, \qquad \tilde{B} := Q^T U_1 B.$$
We observe that $\tilde{u} \in C^2(\tilde{B}; \mathbb{R}^n)$, that $\tilde{U}_1 = Q \in SO(n)$, and that the new matrices $\{Q, \tilde{U}_2, \ldots, \tilde{U}_l\}$ obey (1.5) with $\tilde{\xi}_1 = e_n$. Further, there are $r_0, r_1 > 0$ such that $B(0, r_0) \subset \tilde{B} \subset B(0, r_1)$.

We show that Proposition 3.4 follows from the application of Lemma 3.5 below to the new function $\tilde{u}$. Indeed, we observe that for any $R \in SO(n)$ and any $i$ we have $|\nabla \tilde{u} - R\tilde{U}_i| = |\nabla u U_1^{-1} Q - RU_i U_1^{-1} Q| = |(\nabla u - RU_i)U_1^{-1}| \leq |\nabla u - RU_i|\, |U_1^{-1}|$, and compute
$$\int_{\tilde{B}} \text{dist}(\nabla \tilde{u}, \tilde{K})\, dx \leq \det U_1 \int_B |U_1^{-1}| \text{dist}(\nabla u, K)\, dx \leq c\varepsilon,$$
and analogously, recalling (3.6), (3.9), and $\varepsilon \leq \eta$,
$$\int_{\tilde{B}} \left[\text{dist}(\nabla \tilde{u}, SO(n)) + |\nabla^2 \tilde{u}|\right] dx \leq \det U_1 \int_B \left[|U_1^{-1}|\text{dist}(\nabla u, SO(n)U_1) + |U_1^{-1}|^2 |\nabla^2 u|\right] dx \leq c\eta.$$
Clearly the same estimate holds in $B(0, r_0) \subset \tilde{B}$. Therefore, if $\eta$ is sufficiently small, Lemma 3.5 shows that there is a $\tilde{\rho} > 0$ such that
$$\int_{B(0, \tilde{\rho}r_0)} \text{dist}(\nabla \tilde{u}, SO(n))\, dx \leq c \int_{B(0, r_0)} \text{dist}(\nabla \tilde{u}, \tilde{K})\, dx \leq c \int_B \text{dist}(\nabla u, K)\, dx.$$
Finally, setting $\rho = \tilde{\rho}r_0/r_1$ we obtain $Q^T U_1 B(0, \rho) = \rho\tilde{B} \subset B(0, \tilde{\rho}r_0)$ and
$$\int_{B(0, \rho)} \text{dist}(\nabla u, SO(n)U_1)\, dx \leq c \int_{B(0, \tilde{\rho}r_0)} \text{dist}(\nabla \tilde{u}, SO(n))\, dx \leq c \int_B \text{dist}(\nabla u, K)\, dx$$
which concludes the proof. Notice that all constants depend only on $U_1$, therefore taking the maximum among the $l$ possible majority phases they can be made universal. □

**Lemma 3.5.** Let $U_1 \ldots U_l$ be matrices in $\mathbb{R}^{n \times n}$ with positive determinant, which obey (1.5) with $\xi_1 = e_n$, and assume that $U_1 \in SO(n)$. Let $K$ be defined as in (1.3). Then there are positive numbers $\eta$, $c^*$ and $\rho \in (0,1)$ depending only on $\{U_i\}_i$ and $n$ such that for any $r > 0$, any $u \in C^2(B(0,r); \mathbb{R}^n)$ satisfying

$$\int_{B(0,r)} \mathrm{dist}(\nabla u, SO(n)) + r|\nabla^2 u|\, dx \leq \eta r^n, \tag{3.10}$$

one has

$$\int_{B(0,\rho r)} \mathrm{dist}(\nabla u, SO(n))\, dx \leq c^* \int_{B(0,r)} \mathrm{dist}(\nabla u, K)\, dx.$$

*Proof.* By scaling we can assume $r = 1$. Let $B = B(0,1)$ and

$$\varepsilon := \int_{B(0,1)} \mathrm{dist}(\nabla u, K)\, dx.$$

If $\eta \leq \varepsilon$ there is nothing to prove, hence we can assume $\varepsilon < \eta$. Define a partition $(E_i)_i$ of $B$ by setting

$$E_i := \{x \in B \setminus \cup_{k=1}^{i-1} E_k : \mathrm{dist}(\nabla u, SO(n)U_i) = \mathrm{dist}(\nabla u, K)\}$$

for $i = 1, \ldots, l$ (for $i = 1$ the union is understood to be empty). By the triangular inequality this gives

$$\mathrm{dist}(\nabla u, SO(n)) \leq \mathrm{dist}(\nabla u, K) + \sum_{i=2}^{l} \mathrm{dist}(SO(n), SO(n)U_i)\chi_{E_i}, \tag{3.11}$$

where $\chi_{E_i}$ is the characteristic function of $E_i$. Since the second term is bounded, after redefining $\rho$ it is sufficient to show that there are $\rho > 0$ and $c > 0$ such that

$$\mathcal{L}^n\left(B(0, \frac{\rho}{2}) \setminus E_1\right) \leq c\varepsilon. \tag{3.12}$$

Set $\phi(F) := \mathrm{dist}(F, SO(n))$ for $F \in \mathbb{R}^{n \times n}$. Then (3.10) gives

$$\int_B \phi(\nabla u)\, dx + \int_B |D[\phi(\nabla u)]| \leq c\eta. \tag{3.13}$$

Let $\delta > 0$ be such that

$$|U_i \xi| < 1 - 2\delta \quad \text{for all } \xi \in \mathbb{R}^n \quad \text{such that } |\xi - e_n| < 2\delta \text{ and all } i = 2, \ldots, l. \tag{3.14}$$

Let $c > 0$, $\rho > 0$, $\mu > 0$, $\{a_h\}_h$ and $\bar{a}$ as in Lemma 3.2 ($\mu$ is smaller than $c\delta < 1/2$ and $\rho < 1/4(n+1)$). We define $\tilde{a}_h := a_h - \bar{a}$ so that the symplex with vertices $\{\tilde{a}_h\}$ has barycenter in $0$, and has the same properties as the other one. Possibly reducing $\rho$, we can assume that $B(\tilde{a}_i, 2\rho) \in B$ for all $i$. Fix a small $\theta > 0$. In the following we shall denote be $c_\theta$ the constants which may depend on $\theta$. Let $\eta_\theta$ be such that the condition $\eta < \eta_\theta$ and (3.13) permit to apply Proposition 2.1 to each pair $(\tilde{a}_h, \tilde{a}_k)$, with $\alpha = \rho$ and the chosen $\theta$. Then, for each pair $(h, k)$, $h \neq k$, there is a set $\omega_{hk} \subset B(\tilde{a}_h, \rho) \times B(\tilde{a}_k, \rho)$ such that

$$1 - c_\theta \varepsilon \leq \frac{|u(x) - u(y)|}{|x - y|} \leq 1 + c_\theta \varepsilon \quad \text{for all } (x, y) \in B(\tilde{a}_h, \rho) \times B(\tilde{a}_k, \rho) \setminus \omega_{hk}, \tag{3.15}$$

with $\mathcal{L}^{2n}(\omega_{hk}) \leq \theta$. Summing over the $(n+1)(n+2)/2$ possible pairs, we obtain that the "bad" set

$$\omega := \left\{ (b_1, \ldots b_{n+1}) \in \prod_{h=1}^{n+1} B(\tilde{a}_h, \rho) : \text{ there are } h, k \text{ such that } (b_h, b_k) \in \omega_{hk} \right\} \tag{3.16}$$

satisfies $\mathcal{L}^{n(n+1)}(\omega) \leq c\theta$. Fix now $(b_1, \ldots b_{n+1}) \in \prod_{h=1}^{n+1} B(\tilde{a}_h, \rho) \setminus \omega$. Then (3.15) yields that the $(n+1)n/2$ lengths are preserved under $u$ up to errors of order $\varepsilon$, and by Lemma 3.3 this implies that there is an isometry $x \mapsto I(x) := Qx + b$, $Q \in O(n)$, such that

$$|u(x) - I(x)| \leq c\varepsilon \qquad \text{for } x \in \{b_1, \ldots, b_{n+1}\} \tag{3.17}$$

(in fact, one can show that $Q \in SO(n)$, for most choices; this will however not be needed below). The isometry depends on the choice of the $(b_1, \ldots b_{n+1})$, but the constants here and in all following estimates do not.

For any $p \in B_\rho := B(0, \rho)$, any $h \in \{1, \ldots, n+1\}$, consider the segment $[b_h, p]$ and let $\tau_h := (b_h - p)/|b_h - p|$. By Lemma 3.2 we get $|\tau_h - e_n| \leq \delta$, and, by (3.14), $|U_i \tau_h| < 1 - \delta$ for all $i \geq 2$. Thus we have pointwise

$$|\partial_{\tau_h} u| \leq 1 + \sum_{i=2}^{l} (|U_i \tau_h| - 1)\chi_{E_i} + \text{dist}(\nabla u, K) \leq 1 - \delta \chi_{(B \setminus E_1)} + \text{dist}(\nabla u, K),$$

and integrating over $[b_h, p]$, we get $|u(p) - u(b_h)| \leq |p - b_h| + \sigma_h(p) - \delta_h(p)$, where

$$\sigma_h(p) := \int_{[b_h, p]} \text{dist}(\nabla u, K) \, d\mathcal{H}^1, \qquad \delta_h(p) := \delta \int_{[b_h, p]} \chi_{(B \setminus E_1)} \, d\mathcal{H}^1.$$

Recalling (3.17) we have

$$|u(p) - I(b_h)| \leq |u(p) - u(b_h)| + |u(b_h) - I(b_h)| \leq |p - b_h| + \sigma_h(p) - \delta_h(p) + c\varepsilon.$$

By the second statement in Lemma 3.2, applied with $q = I^{-1}(u(p))$, we have

$$\sum_{h=1}^{n+1} c_h |p - b_h| \leq \sum_{h=1}^{n+1} c_h |I^{-1}(u(p)) - b_h| = \sum_{h=1}^{n+1} c_h |u(p) - I(b_h)|.$$

This implies

$$\sum_{h=1}^{n+1} c_h (\sigma_h(p) - \delta_h(p) + c\varepsilon) \geq 0.$$

Since $1/c \leq c_h \leq c$, with $c$ not depending on $p$ or the $b_h$, we conclude that

$$\sum_{h=1}^{n+1} \delta_h(p) \leq c\varepsilon + c \sum_{h=1}^{n+1} \sigma_h(p).$$

We integrate over all $p \in B_\rho$, and obtain

$$\sum_{h=1}^{n+1} \int_{B_\rho} \delta_h(p) \, dp \leq c\varepsilon + c \sum_{h=1}^{n+1} \int_{B_\rho} \sigma_h(p) \, dp. \tag{3.18}$$

We now show that the $b_h$'s can be chosen so that $\int_{B_\rho} \sigma_h(p) dp \leq c\varepsilon$ for any $h = 1, \ldots, n+1$. Let $f := \text{dist}(\nabla u, K)$ in $B = B(0, 1)$ and 0 elsewhere. Then $\sigma_h(p) = \int_{[b_h, p]} f \, d\mathcal{H}^1 \leq \int_{b_h + \mathbb{R}(p - b_h)} f \, d\mathcal{H}^1$. We integrate over all $p \in B_\rho$, and change to polar coordinates centered in $b_h$ according to $p = b_h + r\nu$, $\nu \in \mathbb{S}^{n-1}$. Since $\frac{1}{n+1} - \rho \leq |b_h| \leq \frac{n}{n+1} + \mu + \rho$ (recall that $\mu < 1/2$ and $\rho < 1/4(n+1)$), we have $B_\rho \subset B(b_h, 3/2) \setminus B(b_h, 1/2(n+1))$, that is the range of $r$ is bounded away from 0 and $\infty$. Further, the integral over the segment $[b_h, b_h + r\nu]$ can be extended to the entire half line $b_h + \nu(0, \infty)$, so that it does not depend on $r$. We can then perform the $r$ integration and change variables again according to $x = b_h + t\nu$, to get

$$\begin{aligned}
\int_{B_\rho} \sigma_h(p) \, dp &\leq \int_{\frac{1}{2(n+1)} \leq r \leq \frac{3}{2}} r^{n-1} \int_{\mathbb{S}^{n-1}} \int_0^\infty f(b_h + t\nu) dt \, d\nu \, dr \\
&= c \int_{\mathbb{R}^n} f(x) \frac{1}{|x - b_h|^{n-1}} \, dx.
\end{aligned}$$

Integrating over all $b_h \in B(\tilde{a}_h, \rho)$ yields

$$\int_{B(\tilde{a}_h,\rho)} \int_{B_\rho} \sigma_h(p)\, dp\, db_h = c \int_{\mathbb{R}^n} f(p) \int_{B(\tilde{a}_h,\rho)} \frac{1}{|p - b_h|^{n-1}}\, db_h\, dp \leq c \int_{\mathbb{R}^n} f(x)\, dx = c\varepsilon.$$

Thus we have that

$$\int_{B_\rho} \sigma_h(p)\, dp \leq c_\theta \varepsilon.$$

for all $b_h \in B(\tilde{a}_h, \rho) \setminus \tilde{\omega}_h$, with $\mathcal{L}^n(\tilde{\omega}_h) \leq \theta$. We observe that the "bad" set $\omega \cup \{(c_1, \ldots, c_{n+1}) : c_h \in \tilde{\omega}_h \text{ for some } h\}$ has measure bounded by $\bar{c}\theta$. Therefore choosing $\theta \leq (\mathcal{L}^n(B_\rho))^{n+1}/2\bar{c}$ the complement is nonempty, and we can choose the $b_h$'s will all desired properties. Therefore (3.18) becomes

$$\sum_{h=1}^{n+1} \int_{B_\rho} \delta_h(p)\, dp \leq c\varepsilon.$$

It remains to show that $\int_{B_\rho} \delta_h dp$ controls the volume of the minority phases (in this final step we only need to consider one of the values of $h$, say, $h = 1$). Changing variables as above, we obtain

$$\int_{B_\rho} \delta_h(p)\, dp = \int_{\frac{1}{2(n+1)} \leq r \leq \frac{3}{2}} r^{n-1} \int_{\mathbb{S}^{n-1}} \chi_{B_\rho}(b_h + r\nu) \int_0^r \chi_{B \setminus E_1}(b_h + t\nu)\, dt\, d\nu\, dr.$$

We bound $r^{n-1}$ by a constant, swap the $t$ with the $r$ integration, and obtain

$$\int_{B_\rho} \delta_h(p)\, dp \geq c \int_{\mathbb{S}^{n-1}} \int_0^1 \chi_{B \setminus E_1}(b_h + t\nu) \int_t^1 \chi_{B_\rho}(b_h + r\nu)\, dr\, dt\, d\nu.$$

If $b_h + t\nu \in B_{\rho/2}$, then the integral in $r$ extends over a length at least $\rho/2$, therefore

$$\int_t^1 \chi_{B_\rho}(b_h + r\nu) dr \geq c\chi_{B_{\rho/2}}(b_h + t\nu).$$

Inserting into the previous estimate and changing variables back we conclude that

$$\begin{aligned}
\int_{B_\rho} \delta_h(p)\, dp &\geq c \int_{\mathbb{S}^{n-1}} \int_0^1 \chi_{B \setminus E_1}(b_h + t\nu) \chi_{B_{\rho/2}}(b_h + t\nu)\, dt\, d\nu \\
&\geq c \int_{B_{\rho/2}} \chi_{B \setminus E_1}(x)\, dx = c\mathcal{L}^n(B_{\rho/2} \setminus E_1).
\end{aligned}$$

This concludes the proof of (3.12) and therefore of Proposition 3.4. □

We finally come to the proof of Theorem 1.1 which is based on Proposition 3.4 and a covering argument of the domain $\Omega$ with suitable balls. The argument given in [9, Pf. of Th. 2.1] can be applied with minimal changes to the $n$-dimensional case, after passing in Proposition 3.4 from balls to cubes. For the convenience of the reader we prefer to report the full argument in the current notation, and in a version which only uses balls. We remark that the exponent of the distance $\text{dist}(\Omega', \partial \Omega)$ in (1.6) comes from (3.6).

*Proof of Theorem 1.1.* Without loss of generality we can assume $SO(n)U_i \cap SO(n)U_j = \emptyset$. We let $\eta > 0$, $\rho \leq 1/2$ and $c^*$ be as in Proposition 3.4 (they all depend only on $K$), and set $d := \text{dist}(\Omega', \partial \Omega)$,

$$c_2 := \frac{1}{4c_* d^n} \mathcal{L}^n\left(B\left(0, \frac{\rho d}{2}\right)\right) \min_{i \neq j} \{\text{dist}(U_i, SO(n)U_j)\}.$$

Assume first that (1.8) does not hold. Then, by (1.6) and the Poincaré inequality,

$$\int_\Omega |\nabla u - F|\, dx \leq c_\Omega \eta_0 d$$

for some $F \in \mathbb{R}^{n\times n}$, and the proof is concluded following the argument used in the proof of Proposition 3.4, case $\eta \leq \varepsilon$.

Assume now that (1.8) holds. Choose $\{p_1,\ldots,p_k\} \in \Omega'$ be such that the closures of the balls $B_h := B(p_h, \rho d/2) \subset \Omega$ cover $\Omega'$, and each point in $\mathbb{R}^n$ belongs to at most a fixed number $M$ (depending on the dimension $n$ and on $\rho$, hence only on $K$) of the larger balls $B'_h := B(p_h, d)$. This can be done using the Besicovitch covering theorem, see, e.g., [8, Lemma 2.2 and 2.3]. Clearly $B'_h \subset \Omega$ for all $h$. Then (1.6) (with $\eta_0 = \eta$) implies that Proposition 3.4 can be applied to each of the $B'_h$; hence for each $h$ there is $J_h \in \{U_i\}_i$ such that

$$\int_{B(p_h, \rho d)} \mathrm{dist}(\nabla u, SO(n)J_h)\, dx \leq c^* \int_{B'_h} \mathrm{dist}(\nabla u, K)\, dx. \qquad (3.19)$$

The key observation is that $J_h$ does not depend on $h$. Since $\Omega'$ is connected it suffices to prove that we can define $J : \Omega' \to K$ so that $J$ is constant in a neighborhood of any $x \in \Omega'$, and $J(x) = J_h$ whenever $x \in B_h \cap \Omega'$. In turn, to do this it suffices to show that for any $x \in \Omega'$, and any $h, k$ with $x \in \overline{B}_h \cap \overline{B}_k$, one has $J_h = J_k$. Since $x \in \overline{B}_h = \overline{B}(p_h, \rho d/2)$, we have $B(x, \rho d/2) \subset B(p_h, \rho d)$, and the same for $k$. Therefore (3.19) implies

$$\begin{aligned}
\mathcal{L}^n(B(x, \rho d/2)) \mathrm{dist}(J_h, SO(n) J_k) &\leq \int_{B(x, \rho d/2)} \mathrm{dist}(\nabla u, SO(n)J_h) + \mathrm{dist}(\nabla u, SO(n)J_k)\, dy \\
&\leq 2c^* \int_{B'_h \cup B'_k} \mathrm{dist}(\nabla u, K)\, dy.
\end{aligned}$$

Estimating the right-hand side with the integral over $\Omega$, and recalling (1.8) and the above choice for $c_2$, shows that $J_h = J_k$.

Finally, the theorem follows by summing over all balls,

$$\begin{aligned}
\int_{\Omega'} \mathrm{dist}(\nabla u, SO(n)J)\, dx &\leq \sum_h \int_{B_h} \mathrm{dist}(\nabla u, SO(n)J)\, dx \\
&\leq c^* \sum_h \int_{B'_h} \mathrm{dist}(\nabla u, K)\, dx \leq M c^* \int_\Omega \mathrm{dist}(\nabla u, K)\, dx
\end{aligned}$$

where we used that any point of $\Omega$ belongs at most to $M$ of the larger balls $B'_h$. $\square$

We finally come back to the separation condition (1.5), and show that it holds in most cases of physical interest. In particular, in materials undergoing solid-solid phase transitions the set $K$ is obtained from a single eigenstrain under the left action of $SO(n)$ and the right action of the point group of the austenitic phase, which is a finite subgroup of $SO(n)$. This structural condition implies (1.5).

**Remark 3.6.** Assume that the set $K$ has the form

$$K = \{Q\bar{U}P : Q \in SO(n),\ P \in \mathcal{P}\} \qquad (3.20)$$

where $\mathcal{P}$ is a finite subgroup of $SO(n)$, and $\bar{U} \in \mathbb{R}^{n\times n}$ a fixed matrix, with $\det \bar{U} > 0$. Then condition (1.5) holds.

To see this, we let $U_1 \ldots U_l$ be such that $K = SO(n)U_1 \cup \cdots \cup SO(n)U_l$, with $SO(n)U_i \cap SO(n)U_j = \emptyset$ for $i \neq j$, and define $\psi_i(\xi) := |U_i \xi|^2$. For every pair $(i,j)$ and almost every $\xi \in \mathbb{S}^{n-1}$ one has $\psi_i(\xi) \neq \psi_j(\xi)$; therefore there is a $\xi \in \mathbb{S}^{n-1}$ such that $\psi_i(\xi) \neq \psi_j(\xi)$ for any $i \neq j$. Let $k \in \{1\ldots l\}$ be such that $\psi_k(\xi) := \max_j \psi_j(\xi)$. Then $|U_k\xi| > |U_h\xi|$ for any $h \neq k$. Pick some $i \in \{1\ldots l\}$, let $P_{ik} \in \mathcal{P}$ be such that $SO(n)U_k = SO(n)U_iP_{ik}$. We claim that $\xi_i = P_{ik}\xi$ satisfies (1.5). Indeed, $|U_i\xi_i| = |U_iP_{ik}\xi| = |U_k\xi|$, whereas for any $j \neq i$ one has $|U_j\xi_i| = |U_jP_{ik}\xi|$. Clearly $U_jP_{ik} \in K$. Since $U_j \notin SO(n)U_i$, we obtain $U_jP_{ik} \notin SO(n)U_iP_{ik} = SO(n)U_k$, hence $U_jP_{ik} \in SO(n)U_h$ for some $h \neq k$, and $|U_kP_{ik}\xi| = |U_h\xi| < |U_k\xi|$. This concludes the proof.

As a special case, this applies to the three-well problem in three dimensions:

**Remark 3.7.** Let $\lambda > 0$, $\lambda \neq 1$, and let $U_1, U_2, U_3$ be the three diagonal matrices with eigenvalues $(1,1,\lambda)$. Then (1.5) holds (this is an immediate consequence of Remark 3.6).

In the case of two matrices the condition (1.5) is even weaker. It holds whenever none of the inequalities $U_1 \leq U_2$ or $U_2 \leq U_1$, in the sense of symmetric matrices, holds; in particular (1.5) holds if the two matrices have the same determinant.

**Remark 3.8.** Let $U_1, U_2$ be such that $U_1 \notin SO(n)U_2$, $\det U_1 = \det U_2 \neq 0$. Then (1.5) holds. Indeed, let $B$ be the unit ball in $\mathbb{R}^n$. If the two ellipsoids $U_1^{-1}B$ and $U_2^{-1}B$ were equal, then $U_1U_2^{-1}$ would be a unit-determinant isometry, i.e., belong to $SO(n)$, against the assumption. The determinant condition shows that the two ellipsoids have the same volume, hence none is a subset of the other. Therefore there are $v \in U_1^{-1}B \setminus U_2^{-1}B$, and $w \in U_2^{-1}B \setminus U_1^{-1}B$. The vectors $v$, $w$ have the stated properties.

A simple example where (1.5) does not hold is provided by the set $K = SO(n) \cup 2SO(n)$.

**Remark 3.9.** If the assumption (1.5) is dropped, then Theorem 1.1 does not hold. Consider for example the set (with $l = 3$ wells) $K = SO(2)\{\text{Id}, \text{Id} + e_1 \otimes e_2, \text{Id} - e_1 \otimes e_2\}$, in two dimensions. Set $\Omega = B(0,2)$, $\Omega' = B(0,1)$. Assume Theorem 1.1 would hold. Pick $\ell \in (0, 1/2)$ so that $20\ell$ is less than the right-hand side of (1.6). Define, for $\varepsilon \in (0, 1/8)$,

$$Q_\varepsilon = \{x \in \mathbb{R}^2 : |x_1| + |x_2/\varepsilon| < \ell\}.$$

We set

$$u_\varepsilon(x) = \begin{cases} x + (\varepsilon\ell - |x_1|\varepsilon - |x_2|)e_1 & \text{if } x \in Q_\varepsilon, \\ x & \text{else.} \end{cases}$$

We observe that $u_\varepsilon \in \text{Lip}(\Omega; \mathbb{R}^2)$, that on $Q_\varepsilon$ we have $\nabla u_\varepsilon = \text{Id} \pm e_1 \otimes e_2 \pm \varepsilon e_1 \otimes e_1$, and $\text{dist}(\nabla u_\varepsilon, K) \leq \varepsilon$ uniformly. We conclude that

$$\int_\Omega \text{dist}(\nabla u_\varepsilon, K)\, dx \leq \varepsilon \mathcal{L}^2(Q_\varepsilon) = 2\varepsilon^2\ell^2,$$

whereas, since $\mathcal{L}^2(\Omega' \setminus Q_\varepsilon) > 2\mathcal{L}^2(Q_\varepsilon)$,

$$\min_{J \in K} \int_{\Omega'} \text{dist}(\nabla u_\varepsilon, SO(2)J)\, dx = \int_{Q_\varepsilon} \text{dist}(\nabla u_\varepsilon, SO(2))\, dx \geq c\mathcal{L}^2(Q_\varepsilon) = 2c\varepsilon\ell^2.$$

Taking $\varepsilon \to 0$ we see that there is no $c_0$ such that (1.7) can hold. The same construction can easily be extended to higher dimension. With a larger set $K$ it is easy to find examples where the right-hand side of (1.7) vanishes, but the left-hand side does not (e.g., $K = SO(2)\nabla u_\varepsilon(\Omega)$).

# Acknowledgements

We thank Bernd Kirchheim for useful discussions. This work was partially supported by the Deutsche Forschungsgemeinschaft through the Schwerpunktprogramm 1095 *Analysis, Modeling and Simulation of Multiscale Problems*.

# References

[1] L. Ambrosio, N. Fusco, and D. Pallara, *Functions of bounded variation and free discontinuity problems*, Mathematical Monographs, Oxford University Press, 2000.

[2] J. M. Ball and R. D. James, *Fine phase mixtures as minimizers of the energy*, Arch. Ration. Mech. Analysis **100** (1987), 13–52.

[3] ______, *Proposed experimental tests of a theory of fine microstructure and the two-well problem*, Phil. Trans. R. Soc. Lond. A **338** (1992), 389–450.

[4] K. Bhattacharya, *Microstructure of martensite: Why it forms and how it gives rise to the shape-memory effect*, Oxford University Press, 2003.

[5] A. Chan and S. Conti, *Energy scaling and domain branching in an $SO(2)$-invariant model*, in preparation.

[6] N. Chaudhuri and S. Müller, *Rigidity estimate for two incompatible wells*, Calc. Var. Partial Differential Equations **19** (2004), 379–390.

[7] S. Conti, I. Fonseca, and G. Leoni, *A $\Gamma$–convergence result for the two-gradient theory of phase transitions*, Comm. Pure Appl. Math. **55** (2002), 857–936.

[8] S. Conti and C. De Lellis, *Sharp upper bounds for a variational problem with singular perturbation*, Math. Ann. **338** (2007), 119–146.

[9] S. Conti and B. Schweizer, *Rigidity and Gamma convergence for solid-solid phase transitions with $SO(2)$-invariance*, Comm. Pure Appl. Math. **59** (2006), 830–868.

[10] ______, *A sharp-interface limit for a two-well problem in geometrically linear elasticity*, Arch. Rat. Mech. Anal. **179** (2006), 413–452.

[11] B. Dacorogna, *Direct methods in the calculus of variations*, Springer-Verlag, New York, 1989.

[12] C. De Lellis and L. Székelyhidi, *Simple proof of two-well rigidity*, C. R. Math. Acad. Sci. Paris **343** (2006), 367–370.

[13] G. Dolzmann, *Variational methods for crystalline microstructure - analysis and computation*, Lecture Notes in Mathematics, no. 1803, Springer-Verlag, 2003.

[14] G. Dolzmann and S. Müller, *Microstructures with finite surface energy: the two-well problem*, Arch. Rat. Mech. Anal. **132** (1995), 101–141.

[15] L. C. Evans and R. F. Gariepy, *Measure theory and fine properties of functions*, Boca Raton CRC Press, 1992.

[16] I. Fonseca and E. Gangbo, *Degree theory in analysis and applications*, Oxford Univ. Press, 1995.

[17] I. Fonseca and C. Mantegazza, *Second order singular perturbation models for phase transitions*, SIAM J. Math. Anal. **31** (2000), no. 5, 1121–1143.

[18] G. Friesecke, R. James, and S. Müller, *A theorem on geometric rigidity and the derivation of nonlinear plate theory from three dimensional elasticity*, Comm. Pure Appl. Math **55** (2002), 1461–1506.

[19] F. W. Gehring, *Rings and quasiconformal mappings in space*, Proc. Nat. Acad. Sci. U.S.A. **47** (1961), 98–105.

[20] F. John, *Rotation and strain*, Comm. Pure Appl. Math. **14** (1961), 391–413.

[21] B. Kirchheim, *Lipschitz minimizers of the 3-well problem having gradients of bounded variation*, Preprint 12, Max Planck Institute for Mathematics in the Sciences, Leipzig, 1998.

[22] R. V. Kohn, *New integral estimates for deformations in terms of their nonlinear strains*, Arch. Rat. Mech. Anal. **78** (1982), 131–172.

[23] A. Lorent, *A two well Liouville theorem*, ESAIM Control Optim. Calc. Var. **11** (2005), 310–356 (electronic).

[24] ______, *The regularisation of the n-well problem by finite elements and by singular perturbation are scaling equivalent*, MPI-MIS Preprint (46/2007).

[25] S. Müller, *Variational models for microstructure and phase transitions*, in: Calculus of variations and geometric evolution problems (F. Bethuel et al., eds.), Springer Lecture Notes in Math. 1713, Springer-Verlag, 1999, pp. 85–210.

[26] M. Pitteri and G. Zanzotto, *Continuum models for phase transitions and twinning in crystals*, CRC/Chapman & Hall, London, 2003.

[27] Yu. G. Reshetnyak, *Liouville's theory on conformal mappings under minimal regularity assumptions*, Sibirskii Math. Z. **8** (1967), 69–85.